# Sufficient Conditions for String Stability

**Iasson Karafyllis[*], Dionysis Theodosis[*] and Markos Papageorgiou[**],[***]**

[*]Dept. of Mathematics, National Technical University of Athens,
Zografou Campus, 15780, Athens, Greece,
emails: iasonkar@central.ntua.gr , dtheodosis@dssl.tuc.gr

[**] Dynamic Systems and Simulation Laboratory,
Technical University of Crete, Chania, 73100, Greece
email: markos@dssl.tuc.gr

[***]Faculty of Maritime and Transportation,
Ningbo University, Ningbo, China

*To Zhong-Ping Jiang on the occasion of his 60th birthday
with gratitude for his highly influential
contribution to control theory*

## Abstract

Zhong-Ping Jiang devoted a large part of his work to the study of the stability properties of interconnected systems. In this short paper we celebrate Zhong-Ping by studying a special class of families of interconnected systems: the so-called strings. We develop trajectory-based and Lyapunov-based tools that allow the verification of string stability to homogeneous bidirectional strings. The obtained results are applied to the problem of cruise controller design for a string of vehicles.

## 1. Introduction

Interconnections are very common both in control systems and in nature. A composite system may be the result of the interconnection of various subsystems through feedback, communication, or shared physical variables. A central question is whether certain desirable properties such as stability, robustness to disturbances, and input-output performance, are preserved or emerge when subsystems are linked via cascade, feedback, or network interconnections.

The study of interconnected systems leverages tools such as Input-to-State Stability (ISS) ( [7], [19], [27], [22], [24]), Input-to-Output Stability (IOS) ([6], [7], [27], [25], [26]), small-gain theorems ([3], [6], [7], [8], [11]), and vector Lyapunov functions ([5], [7], [10], [18]) that exploit the overall system's structure, enabling guarantees that scale with size. These methods can accommodate heterogeneity, nonlinearities, and uncertainties by certifying that interconnection "gains" remain below certain thresholds and that energy-like functions dissipate along trajectories.

Besides Lyapunov stability properties, robustness of interconnected systems is of high importance since it reflects the impact of disturbances and modeling errors on stability properties and performance of the overall system. The main tools that allow the study of robustness irrespective the number of subsystems, are the ISS concept and small-gain theorems, via which, by bounding interconnection gains and disturbance "sizes", one can certify uniform Input-to-State bounds, disturbance attenuation, and margins of uncertainty. The relation between ISS and small-gain plays a key role in addressing the problems of robust stability and stabilization of interconnected systems.

In strings of interconnected systems, robustness focuses on preventing disturbance amplification as the length of the string increases. This property, known as string stability, has been leveraged for modeling, analysis, and controller design in strings of vehicles (platoons). Intuitively, a platoon is said to be string stable if disturbances on the leading vehicle are not amplified when propagating along the string of vehicles (see for instance [1], [2], [13], [28], [30]). This is a very important property and has been often used as a metric for safety, helping vehicles maintain spacing constraints as well as stop-and-go wave dissipation. While several definitions have been used in the literature (see [4] for a review of various definitions), the notion of $L^p$ string stability is very useful since it encompasses upstream disturbance attenuation of the external input of the leading vehicle, as well as perturbations on initial conditions (see [14], [21]).

Zhong-Ping Jiang devoted a large part of his work to the study of the stability properties of interconnected systems. In this short paper, we celebrate Zhong-Ping by studying a special class of families of homogeneous interconnected systems in series: the so-called "strings". More specifically, we provide general sufficient conditions for string stability. Such strings are often encountered in platoons of vehicles with "Follow-the-Leader" or "bidirectional" architecture (see for instance [13], [14], [29], [30]), as well as in the discretization of 1D linear Partial Differential Equations (PDEs; see [9]). The proposed $L^p$ string stability definition covers both bidirectional couplings between subsystems and the one-direction interconnections commonly appearing in vehicle strings (see [14]). We present a unified approach for certifying $L^p$ string stability using both trajectory-based and Lyapunov-based criteria that accommodate upstream and downstream inputs as well as perturbations of initial conditions. By casting the criteria as small-gain inequalities with composite gains strictly below one, our approach applies to general strings of interconnected systems that cover both bidirectional and one-sided interconnections. Therefore, the trajectory-based string stability conditions are reminiscent of the classical small-gain conditions given in Zhong-Ping Jiang's works (see [6]). Finally, the proposed Lyapunov framework is used to study strings of automated vehicles under the effect of decentralized cruise controllers with bidirectional sensing, initially designed in [15], leading to an $L^2$ string-stability guarantee.

The structure of the present paper is as follows. We start with a section that provides the notion of a homogeneous string and the notion of $L^p$ string stability. Section 3 is devoted to the presentation of trajectory-based sufficient conditions for string stability. Lyapunov-based conditions for string stability are discussed in Section 4. While Section 2 and Section 3 contain simple pedagogical examples that illustrate the use of the provided sufficient conditions, we provide in Section 5 an application to a string of vehicles under the effect of cruise controllers. All proofs are given in Section 6 and the concluding remarks of the present work are presented in Section 7.

**Notation and Basic Notions.** Throughout this paper, we adopt the following notation.

* $\mathbb{R}_+ := [0, +\infty)$. For a vector $x \in \mathbb{R}^n$, $|x|$ denotes its Euclidean norm.

* Let $D \subseteq \mathbb{R}^n$ be an open set and let $S \subseteq \mathbb{R}^n$ be a set that satisfies $D \subseteq S \subseteq cl(D)$, where $cl(D)$ is the closure of $D$. By $C^0(S;\Omega)$, we denote the class of continuous functions on $S$, which take values in $\Omega \subseteq \mathbb{R}^m$. By $C^k(S;\Omega)$, where $k \geq 1$ is an integer, we denote the class of functions on $S \subseteq \mathbb{R}^n$, which take values in $\Omega \subseteq \mathbb{R}^m$ and have continuous derivatives of order $k$. In other words, the functions of class $C^k(S;\Omega)$ are the functions which have continuous derivatives of order $k$ in $D = \text{int}(S)$ that can be continued continuously to all points in $\partial D \cap S$. When $\Omega = \mathbb{R}$ then we write $C^0(S)$ or $C^k(S)$. For $V \in C^1(S)$ with $S \subseteq \mathbb{R}^n$, we define
$$\nabla V(x) = \left( \frac{\partial V}{\partial x_1}(x), \ldots, \frac{\partial V}{\partial x_n}(x) \right).$$

* Let $D \subseteq \mathbb{R}^p$ be a non-empty set and let $I \subseteq \mathbb{R}$ be an interval. By $L^\infty(I;D)$ we denote the class of essentially bounded, Lebesgue measurable functions $d: I \to D$. When $D = \mathbb{R}^p$ then we simply write $L^\infty(I)$. For $d \in L^\infty(I;D)$ we define $\|d\|_\infty = \sup_{t \in I}(|d(t)|)$, where $\sup_{t \in I}(|d(t)|)$ is the essential supremum. By $L^\infty_{loc}(\mathbb{R}_+;D)$ we denote the class of essentially bounded, Lebesgue measurable functions $d: \mathbb{R}_+ \to D$ with $d \in L^\infty((0,T);D)$ for all $T > 0$. We define $\mho(D)$ to be the set of functions $u \in L^\infty_{loc}(\mathbb{R}_+;D)$ with the property that for every $T > 0$ there exists a compact set $\Omega_T \subseteq D$ (in the standard topology of $\mathbb{R}^p$) such that $u(t) \in \Omega_T$ for $t \in [0,T]$ a.e. (almost everywhere).

* Let $y \in L^\infty_{loc}(\mathbb{R}_+;\mathbb{R}^n)$ be a given signal. For each $t \geq 0$ and $p \in [1,+\infty)$ we define $\|y\|_{[0,t],p} = \left( \int_0^t |y(s)|^p ds \right)^{1/p}$. We also define $\|y\|_{[0,t],\infty} = \sup_{s \in [0,t]}(|y(s)|)$.

* Let $D \subseteq \mathbb{R}^n$ be a non-empty open set and let a non-empty set $\Omega \subseteq \mathbb{R}^m$ be given. Let $f: D \times \Omega \to \mathbb{R}^n$ be a locally Lipschitz mapping with respect to $x \in D$, i.e., a mapping for which the following property holds: "for every compact set $S \subseteq D \times \Omega$ there exists $L > 0$ such that $|f(x,u) - f(y,u)| \leq L|x-y|$ for all $(x,u) \in S$, $(y,u) \in S$". We say that the control system $\dot{x} = f(x,u)$ is forward complete, if for every $\xi \in D$, $u \in \mho(\Omega)$ the solution $x(t) \in D$ of the initial-value problem $\dot{x} = f(x,u)$ with initial condition $x(0) = \xi$ exists for all $t \geq 0$.

## 2. The Notion of String Stability

Let $D \subseteq \mathbb{R}^k$ be an open set with $0 \in D$ and let two sets $\Omega, S \subseteq D$ with $0 \in S$, $0 \in \Omega$ be given. Let $h: D \to \mathbb{R}^m$ be a continuous function with $h(0) = 0$ and let $f: D \times D \times D \to \mathbb{R}^k$ be a locally Lipschitz mapping. Consider the following "string" of interconnected systems

$$(\Sigma_n): \begin{cases} \dot{x}_i = f(x_{i-1}, x_i, x_{i+1}) \\ y_i = h(x_i) \end{cases}, i = 1, ..., n \text{ with } x_0 = u \in \Omega, x_{n+1} = w \in S. \quad (2.1)$$

Clearly, systems $(\Sigma_n)$ for $n = 1, 2, ...$ constitute a family of control systems parameterized by $n \in \mathbb{N}$ with inputs $(u, w) \in \Omega \times S$. The state of $(\Sigma_n)$ is $x = (x_1, ..., x_n) \in D^n$ and its output is $y = (y_1, ..., y_n) \in \mathbb{R}^{nm}$. We call $(\Sigma_n)$ a "string" of length $n \in \mathbb{N}$. Such a string of interconnected systems is also called a "homogeneous" string because $f$ and $h$ are the same for all component subsystems.

The string (2.1) is called an isotropic homogeneous string when $f(u, x, w) = f(w, x, u)$ for all $x, u, w \in D$ and $\Omega = S$. Otherwise, it is called an anisotropic homogeneous string.

In the literature of strings of vehicles, the reader can frequently encounter the case where $f(x_{i-1}, x_i, x_{i+1}) = f(x_{i-1}, x_i)$ (see [14], [17], [21]). Such a string has a single direction of transmission of the effect of the input $u \in \Omega$ (from 1 to $n$) and is not affected by $w \in S$. That is why this string is called a "one-directional" string. In the general case, called a "bidirectional" string, the string is affected by both inputs $(u, w) \in \Omega \times S$ and the effects of the inputs are transmitted simultaneously to two different directions (from 1 to $n$ and from $n$ to 1).

The inputs $(u, w) \in \Omega \times S$ are external perturbations of the string (2.1) and we are mainly interested in analyzing the effect of the inputs to the state of the string. More specifically, we want to analyze stability properties that are valid *independently of the length of the string*. Such a stability notion is the notion of string stability, which is defined next.

**Definition 1 (String Stability):** *Let $p \in [1, +\infty]$ be given. We say that (2.1) is $L^p$ string stable if there exist non-decreasing functions $a_1, a_2 \in C^0(\mathbb{R}_+; \mathbb{R}_+)$ (independent of $n \in \mathbb{N}$) such that the following properties hold for every $n \in \mathbb{N}$:*

(i) *system (2.1) is forward complete, and*

(ii) *there exists $Q_n \in C^0(D^n; \mathbb{R}_+)$ such that for every $\xi \in D^n$, $u \in \mho(\Omega)$, $w \in \mho(S)$, the solution of the initial value problem (2.1) with $x(0) = \xi$ satisfies the following estimates for all $t \geq 0$ and $i = 1, ..., n$:*

$$\|y_i\|_{[0,t],p} \leq a_1\left(\|\bar{u}\|_{[0,t],p}\right) + a_2\left(\|\bar{w}\|_{[0,t],p}\right) + Q_n(\xi) \quad (2.2)$$

*where $\bar{u} = h(u)$ and $\bar{w} = h(w)$.*

*We say that (2.1) is $L^p$ one-sided string stable when (2.1) with $S = \{0\}$ is $L^p$ string stable.*

**Remarks on Definition 1: (a)** Definition 1 generalizes known string stability notions (see for instance [4], [21]) and can be applied to bidirectional strings. $L^\infty$ string stability notions for bidirectional strings have also appeared in [20], which however mainly focus on the effects of external disturbances/perturbations $d_i$, i.e., $\dot{x}_i = f(x_{i-1}, x_i, x_{i+1}) + d_i$.

**(b)** When $S = \{0\}$ then the input $w \in S$ is absent, i.e., there is no external input from the one side of the string. That is why $L^p$ string stability in this case is called $L^p$ *one-sided string stability*. This is analogous to the definition for $L^p$ string stability presented in [21].

**(c)** String stability essentially requires that the effect of the inputs is not magnified as it is transmitted from one component of the string to another component of the string, uniformly in the string length, for both inputs and initial-condition mismatches.

**(d)** When $p = +\infty$, inequality (2.2) is not an Input-to-Output Stability (IOS; [6], [7], [12], [16], [25], [24]) estimate since IOS would require and estimate of the form

$$|y_i(t)| \le a_1\left(\|\bar{u}\|_{[0,t],p}\right) + a_2\left(\|\bar{w}\|_{[0,t],p}\right) + \beta(|\xi|, t)$$

where $\beta \in KL$ (being continuous, zero at zero and strictly increasing in the first argument, decreasing in the second argument with limit equal to zero as the second argument tends to infinity). On the other hand, in (2.2) the function $Q_n$ is only assumed to be continuous and independent of the time $t$. Thus, IOS is a stronger property than string stability. On the other hand, string stability requires estimates that are independent of the length of the string. This is not the case for IOS.

**(e)** In the absence of inputs ($u \equiv 0$, $w \equiv 0$), inequality (2.2) reduces to $\|y_i\|_{[0,t],p} \le \tilde{Q}_n(\xi)$ for all $t \ge 0$, $i = 1,...,n$, where $\tilde{Q}_n(\xi) = a_1(0) + a_2(0) + Q_n(\xi)$. For $p = +\infty$, this is exactly the Lagrange output stability property, namely, along every solution, each $y_i$ remains bounded by a function of the initial condition only (see [12]).

**(f)** Relation with Input-Output Stability (see [16]). Inequality (2.2) gives an $L^p$ Input-Output stability-type estimate (see Definition 5.1 in [16]). However, string stability (i.e., inequality (2.2)) requires Input-Output stability with the same gain functions for every $n$, namely, a uniform gain which is independent of the dimension of the system, i.e., the length of the string. The main theorems of Section 3, establish that the notion of string stability of Definition 1 imply finite-gain $L^p$ Input-Output stability for the string, with uniform gain for all $n$ (namely, (2.2) holds with $a_1(s) = \tilde{\gamma}s$ and $a_2(s) = \tilde{\delta}s$, where $\tilde{\gamma}, \tilde{\delta}$ are positive constants).

## 3. Trajectory-Based Criteria for String Stability

Since string stability is a stability notion that holds for strings of arbitrary length, it becomes apparent that string stability depends mainly on the characteristics of each component of the string. Our first result shows how one can prove string stability based on the characteristics of the components of the string.

**Theorem 1 (Component-wise trajectory-based sufficient conditions for $L^p$ string stability):** *Suppose that system (2.1) is forward complete for every $n \in \mathbb{N}$. Let $p \in [1, +\infty]$ be given and suppose that there exist $Q \in C^0(D; \mathbb{R}_+)$ and constants $\gamma, \delta \geq 0$ such that for every $\xi \in D$, $u, w \in \mho(D)$, the solution of the initial value problem*

$$\begin{aligned} \dot{x} &= f(u, x, w) \\ y &= h(x) \\ x &\in D \end{aligned} \qquad (3.1)$$

*with $x(0) = \xi$ exists and satisfies the following estimate for all $t \geq 0$:*

$$\|y\|_{[0,t], p} \leq \gamma \|\bar{u}\|_{[0,t], p} + \delta \|\bar{w}\|_{[0,t], p} + Q(\xi) \qquad (3.2)$$

*where $\bar{u} = h(u)$ and $\bar{w} = h(w)$. If $\gamma + \delta < 1$ then (2.1) is $L^p$ string stable and (2.2) holds with $a_1(s) = \dfrac{\gamma s}{1 - \gamma - \delta}$, $a_2(s) = \dfrac{\delta s}{1 - \gamma - \delta}$ for all $s \geq 0$ and $Q_n(\xi_1, ..., \xi_n) = \dfrac{1}{1 - \gamma - \delta} \max_{i=1,...,n} (Q(\xi_i))$ for all $(\xi_1, ..., \xi_n) \in D^n$.*

Theorem 1 generalizes the idea of "having a gain less than one" that is used in several works for bidirectional strings. It is a trajectory-based result, in the sense that (3.2) is an inequality imposed on trajectories of the single-component system (3.1). When the string is essentially one-directional, i.e., when (3.2) holds with $\delta = 0$, we can obtain a slightly different result.

**Theorem 2 (Component-wise trajectory-based sufficient conditions for $L^p$ string stability in essentially one-directional strings):** *Suppose that system (2.1) is forward complete for every $n \in \mathbb{N}$. Let $p \in [1, +\infty]$ be given and suppose that there exist $Q \in C^0(D; \mathbb{R}_+)$ and constants $\gamma, \delta \geq 0$ such that for every $\xi \in D$, $u, w \in \mho(D)$, the solution of the initial value problem (3.1) with $x(0) = \xi$ exists and satisfies estimate (3.2) for all $t \geq 0$. If $\gamma \leq 1$ and $\delta = 0$ then (2.1) is $L^p$ string stable and (2.2) holds with $a_1(s) = \gamma s$ for all $s \geq 0$, $a_2(s) \equiv 0$ and $Q_n(\xi_1, ..., \xi_n) = \sum_{i=1}^{n} Q(\xi_i)$ for all $(\xi_1, ..., \xi_n) \in D^n$.*

It should be noted that simpler versions of Theorem 2 have been known for a long time in the literature. We include it in the present work mainly for completeness purposes.

Theorem 2 studies strings where the input $w \in S$ does not affect the string and the effect of the input is transmitted in a single direction. There is the intermediate case where the input $w \in S$ does not affect the string but the effect of the input can be transmitted to both directions. This is the case of $L^p$ one-sided string stability and the following result deals with this case.

**Theorem 3 (Component-wise trajectory-based sufficient conditions for $L^p$ one-sided string stability):** *Suppose that system (2.1) with $S = \{0\}$ is forward complete for every $n \in \mathbb{N}$. Let $p \in [1,+\infty]$ be given and suppose that there exist $Q \in C^0(D;\mathbb{R}_+)$ and constants $\gamma > 0$, $\delta \geq 0$ such that for every $\xi \in D$, $u, w \in \mho(D)$, the solution of the initial value problem (3.1) with $x(0) = \xi$ exists and satisfies estimate (3.2) for all $t \geq 0$. Define $\theta(\gamma,\delta) = \gamma + \delta$ for $\delta < \gamma$ and $\theta(\gamma,\delta) = 2\sqrt{\gamma\delta}$ for $\delta \geq \gamma$. If $\theta(\gamma,\delta) < 1$ then (2.1) with $S = \{0\}$ is $L^p$ string stable, i.e., (2.1) is $L^p$ one-sided string stable, and (2.2) holds with $a_1(s) = \dfrac{\gamma s}{1-\theta(\gamma,\delta)}$ for all $s \geq 0$, $a_2(s) \equiv 0$ and*

$$Q_n(\xi_1,...,\xi_n) = \frac{\max\left(1,(\gamma^{-1}\delta)^{(n-1)/2}\right)}{1-\theta(\gamma,\delta)} \max_{i=1,...,n}(Q(\xi_i)) \text{ for all } (\xi_1,...,\xi_n) \in D^n.$$

The verification of the trajectory-based estimate (3.2) is crucial for the use of Theorem 1, Theorem 2 and Theorem 3. The following proposition provides the means to verify (3.2) when $p < +\infty$.

**Proposition 1:** *Suppose that system (3.1) is forward complete. Let $p \in [1,+\infty)$ be given and suppose that there exist $V \in C^1(D;\mathbb{R}_+)$ and constants $\gamma, \delta \geq 0$, $c > 0$ such that for every $x, u, w \in D$ the following inequality holds:*

$$\nabla V(x) f(u,x,w) \leq -c|h(x)|^p + \gamma^p c|h(u)|^p + \delta^p c|h(w)|^p \tag{3.3}$$

*Then (3.2) holds with $Q(\xi) = \left(c^{-1}V(\xi)\right)^{1/p}$.*

We conclude this section with a simple example that illustrates the use of all the above results. The example will also be utilized in the following section.

**Example 1:** Consider the string with $D = \mathbb{R}$, arbitrary non-empty sets $\Omega, S \subseteq \mathbb{R}$ with $0 \in S$, $0 \in \Omega$ and $f(u,x,w) = -kx + au + bw$, $h(x) = x$, where $k > 0$, $a, b \in \mathbb{R}$ are constants, i.e., the family of systems parameterized by $n \in \mathbb{N}$:

$$(\Sigma_n): \begin{cases} \dot{x}_i = ax_{i-1} - kx_i + bx_{i+1} \\ y_i = x_i \end{cases}, i = 1,...,n \text{ with } x_0 = u \in \Omega, x_{n+1} = w \in S. \tag{3.4}$$

Linear strings like (3.4) arise frequently when we apply spatial discretization to a 1-D linear Partial Differential Equation (PDE) like

$$\phi_t = \alpha \phi_{xx} + \beta \phi_x + \mu \phi$$

where $\alpha, \beta, \mu \in \mathbb{R}$ are constants (see for example [9]). In such cases, the length $n \in \mathbb{N}$ of the string is large and ideally, we would like to be able to prove properties which are independent of $n \in \mathbb{N}$ because these properties will be inherited by the solution of the PDE.

Since the case $a = b = 0$ is not interesting, we focus on the case where $|a| + |b| > 0$. Without loss of generality, by possible renumbering of the states, we may assume that $a \neq 0$. For system (3.4), forward completeness is not an issue (due to the linearity of (3.4)).

We first consider the case $b = 0$, i.e., the case of having a string which is one-directional. Using Proposition 1 with $p = 2$, $V(x) = \frac{1}{2}x^2$ and the inequality $axu \leq \frac{\varepsilon}{2}|a|x^2 + \frac{1}{2\varepsilon}|a|u^2$, which is valid for arbitrary $u, x \in \mathbb{R}$, $\varepsilon > 0$, we obtain the inequality:

$$\nabla V(x) f(u, x, w) = -kx^2 + axu + bxw$$
$$\leq -\left(k - \frac{\varepsilon}{2}|a|\right)x^2 + \frac{1}{2\varepsilon}|a|u^2$$

We conclude that (3.3) holds with

$$c = k - \frac{\varepsilon}{2}|a| \quad , \quad \gamma = \sqrt{\frac{|a|}{\varepsilon(2k - \varepsilon|a|)}}$$

$$p = 2 \quad , \quad \delta = 0 \quad , \quad Q(\xi) = \frac{|\xi|}{\sqrt{2k - \varepsilon|a|}}$$

for arbitrary $\varepsilon > 0$ with $\varepsilon|a| < 2k$. The selection that minimizes $\gamma$ is $\varepsilon = k/|a|$ and we obtain from Theorem 2 that (3.4) is $L^2$ string stable when

$$b = 0 \text{ and } |a| \leq k. \tag{3.5}$$

Moreover, (2.2) holds with $a_1(s) = |a|s/k$ for all $s \geq 0$, $a_2(s) \equiv 0$ and $Q_n(\xi_1, ..., \xi_n) = \frac{1}{\sqrt{k}} \sum_{i=1}^{n} |\xi_i|$ for all $(\xi_1, ..., \xi_n) \in \mathbb{R}^n$.

We next consider the case $b \neq 0$, i.e., the case of having a string which is bidirectional. Using Proposition 1 with $p = 2$, $V(x) = \frac{1}{2}x^2$ and the inequalities $axu \leq \frac{R}{2}a^2x^2 + \frac{1}{2R}u^2$, $bxw \leq \frac{S}{2}b^2x^2 + \frac{1}{2S}w^2$ that are valid for arbitrary $u, x, w \in \mathbb{R}$, $R, S > 0$, we obtain the inequality:

$$\nabla V(x) f(u, x, w) = -kx^2 + axu + bxw$$
$$\leq -\left(k - \frac{R}{2}a^2 - \frac{S}{2}b^2\right)x^2 + \frac{1}{2R}u^2 + \frac{1}{2S}w^2$$

We conclude that (3.3) holds with

$$c = k - \frac{R}{2}a^2 - \frac{S}{2}b^2 \quad , \quad \gamma = \frac{1}{\sqrt{R(2k - Ra^2 - Sb^2)}}$$

$$p = 2 \quad , \quad \delta = \frac{1}{\sqrt{S(2k - Ra^2 - Sb^2)}} \quad , \quad Q(\xi) = \frac{|\xi|}{\sqrt{2k - Ra^2 - Rb^2}}$$

for arbitrary $R, S > 0$ with $Ra^2 + Sb^2 < 2k$. Since

$$\min\{\gamma + \delta : R > 0, S > 0, Ra^2 + Sb^2 < 2k\}$$
$$= \min\left\{\frac{\sqrt{S} + \sqrt{R}}{\sqrt{RS(2k - Ra^2 - Sb^2)}} : R > 0, S > 0, Ra^2 + Sb^2 < 2k\right\}$$
$$= k^{-1}\left(|b|^{2/3} + |a|^{2/3}\right)^{3/2}$$

with the minimum obtained for $R = \dfrac{k}{|a|^{4/3}\left(|b|^{2/3} + |a|^{2/3}\right)}$, $S = \dfrac{k}{|b|^{4/3}\left(|b|^{2/3} + |a|^{2/3}\right)}$, we conclude from

Theorem 1 that (3.4) is $L^2$ string stable and (2.2) holds with $a_1(s) = \dfrac{|a|^{2/3}\left(|b|^{2/3} + |a|^{2/3}\right)^{1/2}}{k - \left(|b|^{2/3} + |a|^{2/3}\right)^{3/2}} s$,

$a_2(s) = \dfrac{|b|^{2/3}\left(|b|^{2/3} + |a|^{2/3}\right)^{1/2}}{k - \left(|b|^{2/3} + |a|^{2/3}\right)^{3/2}} s$ for all $s \geq 0$ and $Q_n(\xi_1, ..., \xi_n) = \dfrac{k^{1/2} \max\limits_{i=1,...,n}(|\xi_i|)}{k - \left(|b|^{2/3} + |a|^{2/3}\right)^{3/2}}$ for all

$(\xi_1, ..., \xi_n) \in D^n$ when

$$\left(|b|^{2/3} + |a|^{2/3}\right)^{3/2} < k \tag{3.6}$$

Since

$$\min\{\gamma\delta : R \geq S > 0, Ra^2 + Sb^2 < 2k\}$$
$$= \min\left\{\frac{1}{(2k - Ra^2 - Sb^2)\sqrt{RS}} : R \geq S > 0, Ra^2 + Sb^2 < 2k\right\}$$
$$= \begin{cases} 2k^{-2}|ab| & \text{if } |b| \geq |a| \\ k^{-2}(a^2 + b^2) & \text{if } |b| < |a| \end{cases}$$

with the minimum obtained for $R = \dfrac{k}{2a^2}$, $S = \dfrac{k}{2b^2}$ when $|b| \geq a$, we conclude from Theorem 1 that (3.4) with $S = \{0\}$ is $L^2$ string stable, i.e., (3.4) is $L^2$ one-sided string stable, when

$$\left(|b|^{2/3} + |a|^{2/3}\right)^{3/2} < k \text{ and } |b| < |a| \text{ or } 2^{3/2}|ab|^{1/2} < k \text{ and } |b| \geq |a| \tag{3.7}$$

It should be noted that condition (3.7) is less demanding than condition (3.6) when $|b| > |a|$. This is expected, because the conditions for one-sided string stability must be less demanding than the conditions for string stability since one-sided string stability is a weaker stability notion than string stability.

This example will be further studied in the following section. ◁

## 4. Lyapunov-Based Criteria for String Stability

The trajectory-based criteria given in the previous section may be conservative. A reason that explains the conservatism of the trajectory-based criteria is the way that the effect of the input is transmitted: it may happen that the downstream interconnection adds "energy" to a component of the string while the upstream interconnection removes "energy" (or the opposite). In other words, a trajectory-based estimate like (3.2) may ignore the dissipation of the effect of the external inputs.

In many cases, the dissipation of "energy" can be described by appropriate Lyapunov functions. The following result shows how a Lyapunov-based analysis can allow the verification of $L^p$ string stability with $p < +\infty$.

**Theorem 4 (Component-wise Lyapunov-based sufficient conditions for $L^p$ string-stability):** *Suppose that system (2.1) is forward complete for every $n \in \mathbb{N}$. Let $p \in [1, +\infty)$ be given and suppose that there exist $V \in C^1(D; \mathbb{R}_+)$, $g, r \in C^0(D^2; \mathbb{R})$ and constants $\sigma, \omega, \Gamma_1, \Gamma_2 \geq 0$, $c > 0$ such that the following inequalities hold:*

$$\nabla V(x) f(u, x, w) \leq -c |h(x)|^p + g(u, x) + r(x, w), \text{ for all } x, u, w \in D \tag{4.1}$$

$$r(x, w) + g(x, w) \leq c\sigma |h(x)|^p + c\omega |h(w)|^p, \text{ for all } x, w \in D \tag{4.2}$$

$$g(u, x) \leq \Gamma_1 |h(u)|^p + c\omega |h(x)|^p, \text{ for all } x \in D, u \in \Omega \tag{4.3}$$

$$r(x, w) \leq c\sigma |h(x)|^p + \Gamma_2 |h(w)|^p, \text{ for all } x \in D, w \in S \tag{4.4}$$

*If $\sigma + \omega < 1$ then (2.1) is $L^p$ string stable and (2.2) holds with $a_1(s) = \left( \dfrac{\Gamma_1}{c(1-\sigma-\omega)} \right)^{1/p} s$, $a_2(s) = \left( \dfrac{\Gamma_2}{c(1-\sigma-\omega)} \right)^{1/p} s$ for all $s \geq 0$ and $Q_n(\xi_1, ..., \xi_n) = \left( \dfrac{1}{c(1-\sigma-\omega)} \sum_{i=1}^{n} V(\xi_i) \right)^{1/p}$ for all $(\xi_1, ..., \xi_n) \in D^n$. Finally, if*

**(i)** *for every $\rho > 0$ the set $\{ \xi \in D : V(\xi) \leq \rho, |\xi| \leq \rho \}$ is compact (in the standard topology of $\mathbb{R}^k$), and*

**(ii)** *for every compact sets $\bar{\Omega} \subseteq \Omega$, $\bar{S} \subseteq S$ (in the standard topology of $\mathbb{R}^k$) and every $\rho > 0$, the sets $f(\Phi_\rho \times \Phi_\rho \times \Phi_\rho)$, $f(\bar{\Omega} \times \Phi_\rho \times \Phi_\rho)$, $f(\Phi_\rho \times \Phi_\rho \times \bar{S})$ are bounded, where $\Phi_\rho = \{ \xi \in D : V(\xi) \leq \rho \}$,*

*then the assumption that (2.1) is forward complete for every $n \in \mathbb{N}$ is not needed.*

The following result provides Lyapunov-like sufficient conditions for $L^p$ one-sided string stability with $p < +\infty$.

**Theorem 5 (Component-wise Lyapunov-based sufficient conditions for $L^p$ one-sided string-stability):** *Suppose that system (2.1) with $S = \{0\}$ is forward complete for every $n \in \mathbb{N}$. Let $p \in [1, +\infty)$ be given and suppose that there exist $V \in C^1(D; \mathbb{R}_+)$, $g, r \in C^0(D^2; \mathbb{R})$ and constants $\sigma, \omega, \Gamma_1 \geq 0$, $c > 0$, $L \geq 1$ such that (4.1), (4.3) as well as the following inequalities hold:*

$$r(x,w) + Lg(x,w) \leq c\sigma |h(x)|^p + c\omega L |h(w)|^p, \text{ for all } x, w \in D \tag{4.5}$$

$$r(x,0) \leq c\sigma |h(x)|^p, \text{ for all } x \in D \tag{4.6}$$

*If $\sigma + \omega < 1$ then (2.1) with $S = \{0\}$ is $L^p$ string stable and (2.2) holds with*

$$a_1(s) = \left(\frac{\Gamma_1}{c(1-\sigma-\omega)}\right)^{1/p} s \quad \text{for all} \quad s \geq 0, \quad a_2(s) \equiv 0 \quad \text{and}$$

$$Q_n(\xi_1, \ldots, \xi_n) = \left(\frac{1}{c(1-\sigma-\omega)} \sum_{i=1}^n L^{i-1} V(\xi_i)\right)^{1/p} \text{ for all } (\xi_1, \ldots, \xi_n) \in D^n. \text{ Finally, if}$$

**(i)** *for every $\rho > 0$ the set $\{\xi \in D : V(\xi) \leq \rho, |\xi| \leq \rho\}$ is compact (in the standard topology of $\mathbb{R}^k$), and*

**(ii)** *for every compact set $\bar{\Omega} \subseteq \Omega$ (in the standard topology of $\mathbb{R}^k$) and every $\rho > 0$, the sets $f(\Phi_\rho \times \Phi_\rho \times \Phi_\rho)$, $f(\bar{\Omega} \times \Phi_\rho \times \Phi_\rho)$, $f(\Phi_\rho \times \Phi_\rho \times \{0\})$ are bounded, where $\Phi_\rho = \{\xi \in D : V(\xi) \leq \rho\}$,*

*then the assumption that (2.1) with $S = \{0\}$ is forward complete for every $n \in \mathbb{N}$ is not needed.*

There is a big difference between Theorem 4 and Theorem 5. The additional parameter $L \geq 1$ that appears in (4.5) and does not appear in (4.2) allows the relaxation of inequality (4.2). Therefore, as expected, the conditions for one-sided string stability (guaranteed by Theorem 5) are less demanding than the conditions for string stability (guaranteed by Theorem 4).

The following example illustrates the above comment and shows that the Lyapunov-based criteria for string stability provided by Theorem 4 and Theorem 5 are less demanding than the trajectory-based criteria for string stability that were provided in the previous section.

**Back to the Example 1:** Consider again the string with $D = \mathbb{R}$, arbitrary non-empty sets $\Omega, S \subseteq \mathbb{R}$ with $0 \in S$, $0 \in \Omega$ and $f(u,x,w) = -kx + au + bw$, $h(x) = x$, where $k > 0$, $a, b \in \mathbb{R}$ are constants, i.e., the string given by (3.4). Using the function $V(x) = \frac{1}{2}x^2$ and the inequalities

$(a+b)xw \leq \frac{1}{2}|a+b|x^2 + \frac{1}{2}|a+b|w^2$, $aux \leq \frac{a^2}{4k\omega}u^2 + k\omega x^2$, $bxw \leq \frac{b^2}{4k\sigma}w^2 + k\sigma x^2$, that hold for all $\omega, \sigma > 0$, it follows that inequalities (4.1), (4.2), (4.3), (4.4) hold with

$$p = 2 \quad , \quad c = k \quad , \quad g(u,x) = aux$$
$$r(x,w) = bxw \quad , \quad \Gamma_1 = \frac{a^2}{4k\omega} \quad , \quad \Gamma_2 = \frac{b^2}{4k\sigma}$$

and arbitrary $\omega, \sigma > 0$ with $\omega, \sigma \geq \frac{|a+b|}{2k}$. Theorem 4 guarantees that (3.4) is $L^2$ string stable when

$$|a+b| < k \tag{4.7}$$

The reader may notice the difference between (3.6) and (4.7): inequality (4.7) is less demanding than inequality (3.6) since $|a+b| \leq |a| + |b| \leq \left(|b|^{2/3} + |a|^{2/3}\right)^{3/2}$.

Using the function $V(x) = \frac{1}{2}x^2$ and the inequalities $(a+bL)xw \leq \frac{\varepsilon}{2}|aL+b|x^2 + \frac{1}{2\varepsilon}|aL+b|w^2$, $aux \leq \frac{a^2}{4k\omega}u^2 + k\omega x^2$, that hold for all $\varepsilon, \omega > 0$, $L \geq 1$, it follows that inequalities (4.1), (4.3), (4.5), (4.6) hold with

$$p = 2 \quad , \quad c = k \quad , \quad g(u,x) = aux$$
$$r(x,w) = bxw \quad , \quad \Gamma_1 = \frac{a^2}{4k\omega}$$

and arbitrary $L \geq 1$, $\varepsilon, \omega > 0$, $\sigma \geq 0$ with $\sigma \geq \frac{\varepsilon|aL+b|}{2k}$, $\omega \geq \frac{|aL+b|}{2k\varepsilon L}$. Since

$$\min\left\{\left(\varepsilon + \frac{1}{\varepsilon L}\right)|aL+b| : \varepsilon > 0, L \geq 1\right\} = \begin{cases} 0 & \text{if } ab \leq -a^2 \\ 2|a+b| & \text{if } -a^2 < ab < a^2 \\ 4\sqrt{|ab|} & \text{if } ab \geq a^2 \end{cases}$$

it follows from Theorem 5 that system (3.4) with $S = \{0\}$ is $L^p$ string stable, i.e., (3.4) is one-sided string stable when

$$ab \leq -a^2 \text{ or } |a+b| < k \text{ and } -a^2 < ab < a^2 \text{ or } 2\sqrt{|ab|} < k \text{ and } ab \geq a^2 \tag{4.8}$$

The reader can notice the difference between (3.7) and (5.1): inequality (4.8) is less demanding than inequality (3.7) since $|a+b| \leq |a| + |b| \leq \left(|b|^{2/3} + |a|^{2/3}\right)^{3/2}$.

It should be noted that condition (4.8) is less demanding than condition (4.7) when $|b| > |a|$. As noted above this is expected: the conditions for one-sided string stability must be less demanding than the conditions for string stability since one-sided string stability is a weaker stability notion than string stability. ◁

# 5. Application to Cruise Controller Design

In this section we show how we can apply the results of the previous sections to the study of strings of automated vehicles under the effect of decentralized cruise controllers.

A. The Cruise Controller

The longitudinal movement of $n \geq 2$ vehicles on an open road is described by the following set of ODEs:

$$\begin{aligned} \dot{s}_i &= v_{i-1} - v_i, \quad i = 1,...,n \\ \dot{v}_i &= F_i, \quad i = 1,...,n \end{aligned} \tag{5.1}$$

where $s_i$ is the back-to-back distance between vehicle $i$ and its preceding vehicle $i-1$, $v_i$ is the speed of vehicle $i$, and $F_i$ is the acceleration of vehicle $i$. To achieve collision avoidance, we require that for all times $t \geq 0$, the inter-vehicle distances $s_i(t)$, $i=1,...,n$, should be greater than a certain safety distance $L > 0$. In addition, we want to ensure that vehicles never move backwards, i.e., $v_i(t) > 0$, respect the speed limit of the road $v_{max} > 0$, i.e., $v_i(t) < v_{max}$ for all $t \geq 0$, and eventually tend to the same desired speed $v^* \in (0, v_{max})$. The velocity of the leader $v_0 \in (0, v_{max})$ is considered to be an external input. Thus, the state-space of (5.1) is the open set given by

$$D_S = \left\{ \begin{array}{l} (s_1,...,s_n, v_1,...,v_n) \in \mathbb{R}^{2n} : \\ s_i \in (L, +\infty), v_i \in (0, v_{max}), i = 1,...,n \end{array} \right\}. \tag{5.2}$$

The recent papers [13], [15] provided fully decentralized cruise controllers with bidirectional sensing and collision avoidance, by employing artificial potential functions. The cruise controllers in [15] rely only on measurements of distance and speed of both the preceding and following vehicles, when their distance is less than a given interaction distance $\lambda > L$. Let $\Phi \in C^3((L, +\infty); \mathbb{R}_+)$ be a function that satisfies

$$\lim_{x \to L^+} (\Phi(x)) = +\infty, \tag{5.3}$$

$$\Phi(x) = 0, \quad x \geq \lambda. \tag{5.4}$$

$$\begin{aligned} \Phi'(x) &< 0, \quad \text{for } x \in (L, \lambda) \\ \Phi''(x) &> 0, \quad \text{for } x \in (L, \lambda) \end{aligned}. \tag{5.5}$$

We consider the following bidirectional cruise controllers for $i = 1,...,n$

$$F_i = \frac{1}{\beta(v_i, f_i)} \left( v_{max}^2 \frac{Z_i - \mu(v_i - f_i)}{v_i(v_{max} - v_i)} + \Phi'(s_i) - \Phi'(s_{i+1}) \right) \tag{5.6}$$

where $\mu > 0$ is a constant,

$$\beta(v, y) := \frac{v_{max}^3 (v + y) - 2 v_{max}^2 y v}{2 (v_{max} - v)^2 v^2}, \quad v, y \in (0, v_{max}) \tag{5.7}$$

that satisfies $\beta(v, y) > 0$ for all $v, y \in (0, v_{max})$,

$$Z_i = -b'\left(\Phi'(s_{i+1}) - \Phi'(s_i)\right)\left(\Phi''(s_{i+1})(v_i - v_{i+1}) - \Phi''(s_i)(v_{i-1} - v_i)\right) \tag{5.8}$$

$$f_i = v^* - b\left(\Phi'(s_{i+1}) - \Phi'(s_i)\right) \tag{5.9}$$

$v^* \in (0, v_{max})$ is the desired speed of the vehicles, and $b: \mathbb{R} \to (v^* - v_{max}, v^*)$ is a $C^2$ and increasing function satisfying

$$b(0) = 0, \; xb(x) > 0, \; x \neq 0 \text{ and } b'(x) > 0 \text{ for all } x \in \mathbb{R}. \tag{5.10}$$

In the above formulas we use the convention $s_{n+1} = \lambda$ and arbitrary $v_{n+1}$. The design of the cruise controller (4.5), performed in [15], was based on a control Lyapunov methodology.

B. Transformation

The reader can intuitively understand that the closed-loop system (5.1) with (5.6) constitutes a string of vehicles. We next show this fact.

We first apply the following transformation for $i = 1, ..., n$:

$$y_i = \frac{v_i - f_i}{\sqrt{(v_{max} - v_i)v_i}} \tag{5.11}$$

$$z_i = s_i - \lambda$$

The transformation (5.11) maps the set $D_S$ defined by (5.2) onto the set $\mathbb{R}^n \times (L - \lambda, +\infty)^n$. The inverse of transformation (5.11) is given by the following equations for $i = 1, ..., n$:

$$v_i = y_i p\left(y_i, \Phi'(z_{i+1} + \lambda) - \Phi'(z_i + \lambda)\right) + v^* - b\left(\Phi'(z_{i+1} + \lambda) - \Phi'(z_i + \lambda)\right) \tag{5.12}$$

$$s_i = z_i + \lambda$$

where

$$p(y, x) := \frac{y\left(v_{max} - 2v^* + 2b(x)\right) + \sqrt{y^2 v_{max}^2 + 4\left(v_{max} - v^* + b(x)\right)\left(v^* - b(x)\right)}}{2(1 + y^2)} \tag{5.13}$$

The closed-loop system (5.1) with (5.6) under the transformation (5.11) and the following input transformation

$$y_0 = \frac{v_0 - v^* + b\left(\Phi'(s_1)\right)}{\sqrt{(v_{max} - v_0)v_0}} \tag{5.14}$$

is given by the following equations for $i = 1, ..., n$:

$$\dot{z}_i = b\left(\Phi'(z_{i+1} + \lambda) - \Phi'(z_i + \lambda)\right) - b\left(\Phi'(z_i + \lambda) - \Phi'(z_{i-1} + \lambda)\right)$$
$$+ y_{i-1} p\left(y_{i-1}, \Phi'(z_i + \lambda) - \Phi'(z_{i-1} + \lambda)\right) - y_i p\left(y_i, \Phi'(z_{i+1} + \lambda) - \Phi'(z_i + \lambda)\right) \tag{5.15}$$

$$\dot{y}_i = -\mu y_i - \frac{\Phi'(z_{i+1} + \lambda) - \Phi'(z_i + \lambda)}{v_{max}^2} p\left(y_i, \Phi'(z_{i+1} + \lambda) - \Phi'(z_i + \lambda)\right) \tag{5.16}$$

with

$$z_0 \geq 0, y_0 \in \mathbb{R}$$
$$(z_i, y_i) \in (L - \lambda, +\infty) \times \mathbb{R}, i = 1, ..., n \quad (5.17)$$
$$z_{n+1} \geq 0, y_{n+1} \in \mathbb{R}$$

It follows that the system (5.15), (5.16), (5.17) is a string of the form (2.1) with

$$x_i = (z_i, y_i), \; D = (L - \lambda, +\infty) \times \mathbb{R}, \; S = \mathbb{R}_+ \times \mathbb{R}, \; \Omega \subseteq \mathbb{R}_+ \times \mathbb{R} \quad (5.18)$$

$$f(u, x, w) = \begin{bmatrix} \kappa\left(u_y, \Phi'(z+\lambda) - \Phi'(u_z+\lambda)\right) - \kappa\left(y, \Phi'(w_z+\lambda) - \Phi'(z+\lambda)\right) \\ -\mu y - \dfrac{\Phi'(w_z+\lambda) - \Phi'(z+\lambda)}{v_{\max}^2} p\left(y, \Phi'(w_z+\lambda) - \Phi'(z+\lambda)\right) \end{bmatrix},$$

$$\text{for } u = (u_z, u_y), x = (z, y), w = (w_z, w_y) \quad (5.19)$$

where

$$\kappa(y, x) := yp(y, x) - b(x). \quad (5.20)$$

C. String Stability Analysis

We next proceed to the study of the string stability properties of the string (2.1), (5.18), (5.19). To this purpose we employ the function

$$V(z, y) = \frac{v_{\max}^2}{2} y^2 + \Phi(\lambda + z), \text{ for } x = (z, y) \in D = (L - \lambda, +\infty) \times \mathbb{R} \quad (5.21)$$

Applying Theorem 4 with $V$ defined as above, we obtain the following proposition.

**Proposition 6:** *There exist constants $K, \Lambda > 0$ such that (2.1), (5.18), (5.19) with $\Omega = \mathbb{R}_+ \times [-\Lambda, \Lambda]$ is $L^2$ string stable and (2.2) holds with $a_1(s) = Ks$ for all $s \geq 0$, $a_2(s) \equiv 0$ and $Q_n(\xi_1, ..., \xi_n) = \left(\dfrac{1}{\mu v_{\max}^2} \sum_{i=1}^{n} V(\xi_i)\right)^{1/2}$ for all $(\xi_1, ..., \xi_n) \in D^n$.*

There is a big difference between Proposition 6 and Theorem 6 in [15]. Theorem 6 in [15] guarantees an estimate like (2.2) for a restricted class of external inputs, as the allowable size of external inputs depends both on the initial condition as well as on the length of the string. On the other hand, Proposition 6 is a stronger result and guarantees that the allowable size of external inputs is a constant $\Lambda > 0$, independent of the initial condition and the length of the string. The proof of Proposition 6 provides the conservative estimate $\Lambda = \dfrac{2}{5 v_{\max}} \lim_{s \to -\infty} (-b(s))$ for the constant $\Lambda > 0$. However, this estimate of $\Lambda$ is conservative (due to the conservative nature of the inequalities used in the proof of Proposition 6) and external inputs of bigger size can be allowed (see also the simulations in [15]).

## 6. Proofs

We start with the proofs of Theorem 1 and Theorem 2.

**Proof of Theorem 1:** Let arbitrary $\xi = (\xi_1,...,\xi_n) \in D^n$, $u \in \mho(\Omega)$, $w \in \mho(S)$ be given and consider the unique solution $x(t)$ of the initial value problem (2.1) with $x(0) = \xi$. By virtue of (3.2) we obtain the following inequalities for all $t \geq 0$:

$$\|y_1\|_{[0,t],p} \leq \gamma \|\bar{u}\|_{[0,t],p} + \delta \|y_2\|_{[0,t],p} + Q(\xi_1) \tag{6.1}$$

$$\|y_i\|_{[0,t],p} \leq \gamma \|y_{i-1}\|_{[0,t],p} + \delta \|y_{i+1}\|_{[0,t],p} + Q(\xi_i), \text{ for } i = 2,...,n-1 \tag{6.2}$$

$$\|y_n\|_{[0,t],p} \leq \gamma \|y_{n-1}\|_{[0,t],p} + \delta \|\bar{w}\|_{[0,t],p} + Q(\xi_n) \tag{6.3}$$

where $\bar{u} = h(u)$ and $\bar{w} = h(w)$.

Define for $t \geq 0$:

$$Y(t) := \max_{i=1,...,n} \left( \|y_i\|_{[0,t],p} \right) \tag{6.4}$$

Using (6.1), (6.2), (6.3) and (6.4) we obtain for $i = 1,...,n$ and $t \geq 0$:

$$\|y_i\|_{[0,t],p} \leq (\gamma + \delta) Y(t) + \max_{j=1,...,n} \left( Q(\xi_j) \right) + \max \left( \gamma \|\bar{u}\|_{[0,t],p}, \delta \|\bar{w}\|_{[0,t],p} \right) \tag{6.5}$$

Since $\gamma + \delta < 1$, definition (6.4) and estimates (6.5) imply the following estimate for all $t \geq 0$:

$$\begin{aligned}
\max_{i=1,...,n} \left( \|y_i\|_{[0,t],p} \right) &\leq (1-\gamma-\delta)^{-1} \max_{i=1,...,n} \left( Q(\xi_i) \right) \\
&+ (1-\gamma-\delta)^{-1} \gamma \|\bar{u}\|_{[0,t],p} + (1-\gamma-\delta)^{-1} \delta \|\bar{w}\|_{[0,t],p}
\end{aligned} \tag{6.6}$$

Estimates (2.2) with $a_1(s) = \dfrac{\gamma s}{1-\gamma-\delta}$, $a_2(s) = \dfrac{\delta s}{1-\gamma-\delta}$ for $s \geq 0$ and $Q_n(\xi_1,...,\xi_n) = \dfrac{1}{1-\gamma-\delta} \max_{i=1,...,n} \left( Q(\xi_i) \right)$ for $\xi = (\xi_1,...,\xi_n) \in D^n$ are direct consequences of estimate (6.6). The proof is complete. ◁

**Proof of Theorem 2:** Let arbitrary $\xi = (\xi_1,...,\xi_n) \in D^n$, $u \in \mho(\Omega)$, $w \in \mho(S)$ be given and consider the unique solution $x(t)$ of the initial value problem (2.1) with $x(0) = \xi$. By virtue of (3.2) with $\delta = 0$ we obtain the following inequalities for all $t \geq 0$:

$$\|y_1\|_{[0,t],p} \leq \gamma \|\bar{u}\|_{[0,t],p} + Q(\xi_1) \tag{6.7}$$

$$\|y_i\|_{[0,t],p} \leq \gamma \|y_{i-1}\|_{[0,t],p} + Q(\xi_i), \text{ for } i = 2,...,n \tag{6.8}$$

where $\bar{u} = h(u)$. Since $\gamma \leq 1$ we obtain from (6.7) and (6.8) for $i = 1,...,n$ and $t \geq 0$:

$$\|y_i\|_{[0,t],p} \leq \gamma^i \|\bar{u}\|_{[0,t],p} + \sum_{i=1}^{i} Q(\xi_j) \tag{6.9}$$

Estimates (2.2) with $a_1(s) = \gamma s$ for $s \geq 0$, $a_2(s) \equiv 0$ and $Q_n(\xi_1,...,\xi_n) = \sum_{i=1}^{n} Q(\xi_i)$ for $(\xi_1,...,\xi_n) \in D^n$ are direct consequences of estimates (6.9) and the fact that $\gamma \leq 1$ (which implies that $\gamma^i \leq \gamma$ for $i = 1,...,n$). The proof is complete. ◁

We continue with the proofs of Theorem 3 and Proposition 1.

**Proof of Theorem 3:** Let arbitrary $\xi = (\xi_1,...,\xi_n) \in D^n$, $u \in \mho(\Omega)$, $w \in \mho(S)$ be given and consider the unique solution $x(t)$ of the initial value problem (2.1) with $x(0) = \xi$. By virtue of (3.2) and the facts that $h(0) = 0$, $S = \{0\}$, we obtain the following inequalities for all $t \geq 0$:

$$\|y_1\|_{[0,t],p} \leq \gamma \|\bar{u}\|_{[0,t],p} + \delta \|y_2\|_{[0,t],p} + Q(\xi_1) \tag{6.10}$$

$$\|y_i\|_{[0,t],p} \leq \gamma \|y_{i-1}\|_{[0,t],p} + \delta \|y_{i+1}\|_{[0,t],p} + Q(\xi_i), \text{ for } i = 2,...,n-1 \tag{6.11}$$

$$\|y_n\|_{[0,t],p} \leq \gamma \|y_{n-1}\|_{[0,t],p} + Q(\xi_n) \tag{6.12}$$

where $\bar{u} = h(u)$.

Define $L = 1$ for $\delta < \gamma$ and $L = \sqrt{\gamma^{-1}\delta}$ for $\delta \geq \gamma$. Notice that $\gamma L + \frac{\delta}{L} = \theta(\gamma,\delta) = \gamma + \delta$ for $\delta < \gamma$ and $\gamma L + \frac{\delta}{L} = \theta(\gamma,\delta) = 2\sqrt{\gamma\delta}$ for $\delta \geq \gamma$. Moreover, notice that $L \geq 1$. Define for $t \geq 0$:

$$Y(t) := \max_{i=1,...,n} \left( L^i \|y_i\|_{[0,t],p} \right) \tag{6.13}$$

Using (6.10), (6.11), (6.12) and (6.13) we obtain for $i = 1,...,n$ and $t \geq 0$:

$$L^i \|y_i\|_{[0,t],p} \leq \left( \gamma L + \frac{\delta}{L} \right) Y(t) + \max_{j=1,...,n} \left( L^j Q(\xi_j) \right) + \gamma L \|\bar{u}\|_{[0,t],p} \tag{6.14}$$

Exploiting (6.14) and the fact that $\gamma L + \frac{\delta}{L} = \theta(\gamma,\delta) < 1$ we obtain for $t \geq 0$:

$$\max_{i=1,...,n} \left( L^i \|y_i\|_{[0,t],p} \right) \leq \left(1 - \theta(\gamma,\delta)\right)^{-1} \max_{i=1,...,n} \left( L^i Q(\xi_i) \right) + \left(1 - \theta(\gamma,\delta)\right)^{-1} \gamma L \|\bar{u}\|_{[0,t],p} \tag{6.15}$$

It follows from (6.15) and the fact that $L \geq 1$ that the following estimates hold for $i = 1,...,n$ and $t \geq 0$:

$$\|y_i\|_{[0,t],p} \leq (1-\theta(\gamma,\delta))^{-1} L^{n-1} \max_{j=1,\ldots,n}(Q(\xi_j)) + (1-\theta(\gamma,\delta))^{-1} \gamma \|\bar{u}\|_{[0,t],p} \qquad (6.16)$$

Estimates (2.2) with $a_1(s) = \dfrac{\gamma s}{1-\theta(\gamma,\delta)}$ for $s \geq 0$, $a_2(s) \equiv 0$ and

$Q_n(\xi_1,\ldots,\xi_n) = \dfrac{\max\left(1,(\gamma^{-1}\delta)^{(n-1)/2}\right)}{1-\theta(\gamma,\delta)} \max_{i=1,\ldots,n}(Q(\xi_i))$ for all $(\xi_1,\ldots,\xi_n) \in D^n$ are direct consequences of estimates (6.16) and the fact that $L=1$ for $\delta < \gamma$ and $L = \sqrt{\gamma^{-1}\delta}$ for $\delta \geq \gamma$. The proof is complete. ◁

**Proof of Proposition 1:** Let arbitrary $\xi \in D$, $u, w \in \mho(D)$ and consider the solution of the initial value problem (3.1) with $x(0) = \xi$. Using (3.3) we obtain the following differential inequality for $t \geq 0$ a.e.:

$$\dot{V}(x(t)) \leq -c|h(x(t))|^p + \gamma^p c|h(u(t))|^p + \delta^p c|h(w(t))|^p \qquad (6.17)$$

Integrating (6.17) and using the equations $y = h(x)$, $\bar{u} = h(u)$ and $\bar{w} = h(w)$ we get for all $t \geq 0$:

$$V(x(t)) + c\int_0^t |y(s)|^p\, ds \leq V(x(0)) + \gamma^p c \int_0^t |\bar{u}(s)|^p\, ds + \delta^p c \int_0^t |\bar{w}(s)|^p\, ds \qquad (6.18)$$

Since $V$ is non-negative and $x(0) = \xi$, we get from (6.18) for all $t \geq 0$:

$$\|y\|^p_{[0,t],p} \leq c^{-1} V(\xi) + \gamma^p \|\bar{u}\|^p_{[0,t],p} + \delta^p \|\bar{w}\|^p_{[0,t],p} \qquad (6.19)$$

Estimate (3.2) with $Q(\xi) = \left(c^{-1} V(\xi)\right)^{1/p}$ is a direct consequence of (6.19) and the fact that $(a+b)^{1/p} \leq a^{1/p} + b^{1/p}$ for all $a, b \geq 0$. The proof is complete. ◁

We next provide the proofs of Theorem 4 and Theorem 5.

**Proof of Theorem 4:** Let arbitrary $\xi = (\xi_1,\ldots,\xi_n) \in D^n$, $u \in \mho(\Omega)$, $w \in \mho(S)$ be given and consider the unique solution $x(t)$ of the initial value problem (2.1) with $x(0) = \xi$. Then there exists $t_{\max} > 0$ (the maximal existence time of the solution) such that the solution $x(t)$ is defined for all $t \in [0, t_{\max})$. Moreover, if $t_{\max} < +\infty$ then for every non-empty compact set $U \subseteq D^n$ there exists $t \in [0, t_{\max})$ with $x(t) \notin U$ (a consequence of Proposition C.3.6 and Lemma 2.6.2 in [23]; notice that while in Lemma 2.6.2 $f$ is considered $C^1$, Lemma 2.6.2 in [23] can be directly extended to the case where $f$ is locally Lipschitz). Define:

$$W(x) = \sum_{i=1}^n V(x_i) \qquad (6.20)$$

Using (2.1) and (6.20) we get for $t \in [0, t_{\max})$ a.e.:

$$\frac{d}{dt}(W(x(t))) = \sum_{i=1}^{n} \nabla V(x_i(t)) f(x_{i-1}(t), x_i(t), x_{i+1}(t)) \tag{6.21}$$

Exploiting (4.1) and (6.21) we get for $t \in [0, t_{\max})$ a.e.:

$$\begin{aligned}
\frac{d}{dt}(W(x(t))) &\leq -c\sum_{i=1}^{n}|h(x_i(t))|^p + \sum_{i=1}^{n} g(x_{i-1}(t), x_i(t)) + \sum_{i=1}^{n} r(x_i(t), x_{i+1}(t)) \\
&\leq -c\sum_{i=1}^{n}|h(x_i(t))|^p + \sum_{i=0}^{n-1} g(x_i(t), x_{i+1}(t)) + \sum_{i=1}^{n-1} r(x_i(t), x_{i+1}(t)) + r(x_n(t), w(t)) \\
&\leq -c\sum_{i=1}^{n}|h(x_i(t))|^p + \sum_{i=1}^{n-1}\bigl(g(x_i(t), x_{i+1}(t)) + r(x_i(t), x_{i+1}(t))\bigr) \\
&\quad + r(x_n(t), w(t)) + g(u(t), x_1(t))
\end{aligned} \tag{6.22}$$

Using (4.2), we obtain from (6.22) for $t \in [0, t_{\max})$ a.e.:

$$\begin{aligned}
\frac{d}{dt}(W(x(t))) &\leq -c\sum_{i=1}^{n}|h(x_i(t))|^p + c\sum_{i=1}^{n-1}\bigl(\sigma|h(x_i(t))|^p + \omega|h(x_{i+1}(t))|^p\bigr) \\
&\quad + g(u(t), x_1(t)) + r(x_n(t), w(t)) \\
&\leq -c(1-\sigma-\omega)\sum_{i=1}^{n}|h(x_i(t))|^p + g(u(t), x_1(t)) \\
&\quad + r(x_n(t), w(t)) - c\sigma|h(x_n(t))|^p - c\omega|h(x_1(t))|^p
\end{aligned} \tag{6.23}$$

Combining (4.3), (4.4) and (6.23) we get for $t \in [0, t_{\max})$ a.e.:

$$\frac{d}{dt}(W(x(t))) \leq -c(1-\sigma-\omega)\sum_{i=1}^{n}|h(x_i(t))|^p + \Gamma_1|h(u(t))|^p + \Gamma_2|h(w(t))|^p$$

Integrating the above differential inequality and using the equations $y = h(x)$, $\bar{u} = h(u)$ and $\bar{w} = h(w)$ we get for all $t \in [0, t_{\max})$:

$$W(x(t)) + c(1-\sigma-\omega)\sum_{i=1}^{n}\int_0^t |y_i(s)|^p\,ds \leq W(x(0)) + \Gamma_1 \int_0^t |\bar{u}(s)|^p\,ds + \Gamma_2 \int_0^t |\bar{w}(s)|^p\,ds \tag{6.24}$$

The proof now follows two different methodologies, depending on the assumptions in Theorem 4.

1) First, we assume that system (2.1) is forward complete for every $n \in \mathbb{N}$. Since $W$ is non-negative and $x(0) = \xi$, we get from (6.24) for all $t \geq 0$:

$$\sum_{i=1}^{n}\|y_i\|_{[0,t],p}^p \leq \frac{1}{c(1-\sigma-\omega)}\Bigl(W(\xi) + \Gamma_1 \|\bar{u}\|_{[0,t],p}^p + \Gamma_2 \|\bar{w}\|_{[0,t],p}^p\Bigr) \tag{6.25}$$

The fact that estimates (2.2) hold with $a_1(s) = \left( \frac{\Gamma_1}{c(1-\sigma-\omega)} \right)^{1/p} s$, $a_2(s) = \left( \frac{\Gamma_2}{c(1-\sigma-\omega)} \right)^{1/p} s$ for all $s \geq 0$ and $Q_n(\xi_1,...,\xi_n) = \left( \frac{1}{c(1-\sigma-\omega)} \sum_{i=1}^{n} V(\xi_i) \right)^{1/p}$ for all $(\xi_1,...,\xi_n) \in D^n$ is a direct consequence of (6.25), definition (6.20) and the fact that $(a+b)^{1/p} \leq a^{1/p} + b^{1/p}$ for all $a,b \geq 0$.

2) We do not assume that system (2.1) is forward complete for every $n \in \mathbb{N}$. Instead, we assume that properties (i), (ii) in the statement of Theorem 4 hold. In this case, first we show that system (2.1) is forward complete for every $n \in \mathbb{N}$ and then we establish estimates (2.2) as in the previous case.

The proof is made by contradiction. Suppose that $t_{\max} < +\infty$ for certain $\xi = (\xi_1,...,\xi_n) \in D^n$, $u \in \mho(\Omega)$, $w \in \mho(S)$. Then it follows from (6.24), the fact that $x(0) = \xi$ and definition (6.20) that $x_i(t) \in \Phi_\rho = \{ \zeta \in D : V(\zeta) \leq \rho \}$ for all $i = 1,...,n$ and $t \in [0, t_{\max})$ with

$$\rho = W(\xi) + \Gamma_1 \int_0^{t_{\max}} |\bar{u}(s)|^p \, ds + \Gamma_2 \int_0^{t_{\max}} |\bar{w}(s)|^p \, ds.$$

Since $u \in \mho(\Omega)$, $w \in \mho(S)$, there exist compact sets $\bar{\Omega} \subseteq \Omega$, $\bar{S} \subseteq S$ (in the standard topology of $\mathbb{R}^k$) such that $u(t) \in \bar{\Omega}$, $w(t) \in \bar{S}$ for $t \in [0, t_{\max}]$ a.e.. It follows from (2.1) that $\dot{x}_i(t) \in f(\Phi_\rho \times \Phi_\rho \times \Phi_\rho)$ for $i = 2,...,n-1$, $\dot{x}_1(t) \in f(\bar{\Omega} \times \Phi_\rho \times \Phi_\rho)$ and $\dot{x}_n(t) \in f(\Phi_\rho \times \Phi_\rho \times \bar{S})$ for $t \in [0, t_{\max})$ a.e..

Since the sets $f(\Phi_\rho \times \Phi_\rho \times \Phi_\rho)$, $f(\bar{\Omega} \times \Phi_\rho \times \Phi_\rho)$, $f(\Phi_\rho \times \Phi_\rho \times \bar{S})$ are bounded (property (i)), there exists $M > 0$ such that $|\dot{x}(t)| \leq M$ for $t \in [0, t_{\max})$ a.e.. It follows (using the fact that $x(0) = \xi$ and the triangle inequality) that $|x(t)| \leq |\xi| + M t_{\max}$. Therefore, $x_i(t) \in \{ \xi \in D : V(\xi) \leq \bar{\rho}, |\xi| \leq \bar{\rho} \}$ with $\bar{\rho} = |\xi| + M t_{\max} + W(\xi) + \Gamma_1 \int_0^{t_{\max}} |\bar{u}(s)|^p \, ds + \Gamma_2 \int_0^{t_{\max}} |\bar{w}(s)|^p \, ds$ for all $i = 1,...,n$ and $t \in [0, t_{\max})$. Since the set $\{ \xi \in D : V(\xi) \leq \bar{\rho}, |\xi| \leq \bar{\rho} \}$ is compact (property (ii)), it follows that the set $U = \{ \xi \in D : V(\xi) \leq \bar{\rho}, |\xi| \leq \bar{\rho} \}^n$ is compact and $x(t) \in U$ for all $t \in [0, t_{\max})$. This is a contradiction with the fact that $t_{\max} < +\infty$, because a finite escape time would imply that for every non-empty compact set $U \subseteq D^n$ there exists $t \in [0, t_{\max})$ with $x(t) \notin U$ (a consequence of Proposition C.3.6 and Lemma 2.6.2 in [23].).
The proof is complete. ◁

**Proof of Theorem 5:** Let arbitrary $\xi = (\xi_1,...,\xi_n) \in D^n$, $u \in \mho(\Omega)$, $w \in \mho(S)$ be given and consider the unique solution $x(t)$ of the initial value problem (2.1) with $x(0) = \xi$. Then there exists $t_{\max} > 0$ (the maximal existence time of the solution) such that the solution $x(t)$ is defined for all $t \in [0, t_{max})$. Moreover, if $t_{max} < +\infty$ then for every non-empty compact set $U \subseteq D^n$ there exists $t \in [0, t_{max})$ with $x(t) \notin U$ (a consequence of Proposition C.3.6 and Lemma 2.6.2 in [23]). Define:

$$W(x) = \sum_{i=1}^{n} L^i V(x_i) \tag{6.26}$$

Using (2.1) and (6.26) we get for $t \in [0, t_{max})$ a.e.:

$$\frac{d}{dt}(W(x(t))) = \sum_{i=1}^{n} L^i \nabla V(x_i(t)) f(x_{i-1}(t), x_i(t), x_{i+1}(t)) \tag{6.27}$$

Exploiting (4.1) and (6.27) we get for $t \in [0, t_{max})$ a.e.:

$$\begin{aligned}\frac{d}{dt}(W(x(t))) &\leq -c \sum_{i=1}^{n} L^i |h(x_i(t))|^p + \sum_{i=1}^{n} L^i \left(g(x_{i-1}(t), x_i(t)) + r(x_i(t), x_{i+1}(t))\right) \\ &\leq -c \sum_{i=1}^{n} L^i |h(x_i(t))|^p + \sum_{i=0}^{n-1} L^{i+1} g(x_i(t), x_{i+1}(t)) + \sum_{i=1}^{n-1} L^i r(x_i(t), x_{i+1}(t)) + L^n r(x_n(t), 0) \\ &\leq -c \sum_{i=1}^{n} L^i |h(x_i(t))|^p + \sum_{i=1}^{n-1} L^i \left(Lg(x_i(t), x_{i+1}(t)) + r(x_i(t), x_{i+1}(t))\right) \\ &\quad + L^n r(x_n(t), 0) + Lg(u(t), x_1(t))\end{aligned} \tag{6.28}$$

Using (4.5), we obtain from (6.28) for $t \in [0, t_{max})$ a.e.:

$$\begin{aligned}\frac{d}{dt}(W(x(t))) &\leq -c \sum_{i=1}^{n} L^i |h(x_i(t))|^p + c \sum_{i=1}^{n-1} L^i \left(\sigma |h(x_i(t))|^p + \omega L |h(x_{i+1}(t))|^p\right) \\ &\quad + Lg(u(t), x_1(t)) + L^n r(x_n(t), 0) \\ &\leq -c(1 - \sigma - \omega) \sum_{i=1}^{n} L^i |h(x_i(t))|^p + Lg(u(t), x_1(t)) \\ &\quad + L^n r(x_n(t), 0) - c\sigma L^n |h(x_n(t))|^p - c\omega L |h(x_1(t))|^p\end{aligned} \tag{6.29}$$

Combining (4.3), (4.6) and (6.29) we get for $t \in [0, t_{max})$ a.e.:

$$\frac{d}{dt}(W(x(t))) \leq -c(1 - \sigma - \omega) \sum_{i=1}^{n} L^i |h(x_i(t))|^p + \Gamma_1 L |h(u(t))|^p \tag{6.30}$$

Integrating (6.30) and using the equations $y = h(x)$ and $\bar{u} = h(u)$ we get for all $t \in [0, t_{max})$:

$$W(x(t)) + c(1 - \sigma - \omega) \sum_{i=1}^{n} L^i \int_0^t |y_i(s)|^p \, ds \leq W(x(0)) + \Gamma_1 L \int_0^t |\bar{u}(s)|^p \, ds \tag{6.31}$$

The proof now follows two different methodologies, depending on the assumptions in Theorem 4.

1) First, we assume that system (2.1) with $S = \{0\}$ is forward complete for every $n \in \mathbb{N}$. Since $W$ is non-negative and $x(0) = \xi$, we get from (6.31) for all $t \geq 0$ and $i = 1, \ldots, n$:

$$\|y_i\|_{[0,t],p}^p \leq \frac{1}{c(1-\sigma-\omega)}\left(L^{-i}W(\xi) + \Gamma_1 L^{1-i}\|\bar{u}\|_{[0,t],p}^p\right) \qquad (6.32)$$

The fact that estimates (2.2) hold with $a_1(s) = \left(\frac{\Gamma_1}{c(1-\sigma-\omega)}\right)^{1/p} s$ for all $s \geq 0$, $a_2(s) \equiv 0$ and

$Q_n(\xi_1, \ldots, \xi_n) = \left(\frac{1}{c(1-\sigma-\omega)}\sum_{i=1}^{n} L^{i-1}V(\xi_i)\right)^{1/p}$ for all $(\xi_1, \ldots, \xi_n) \in D^n$ is a direct consequence of

(6.32), definition (6.26), the fact that $L \geq 1$ and the fact that $(a+b)^{1/p} \leq a^{1/p} + b^{1/p}$ for all $a, b \geq 0$.

2) We do not assume that system (2.1) is forward complete for every $n \in \mathbb{N}$. Instead, we assume that properties (i), (ii) in the statement of Theorem 4 hold. In this case, first we show that system (2.1) is forward complete for every $n \in \mathbb{N}$ and then we establish estimates (2.2) as in the previous case. This is achieved by following exactly the same procedure as in the proof of Theorem 4.
The proof is complete. ◁

Finally, we end this section with the proof of Proposition 6.

**Proof of Proposition 6:** Using (5.19), (5.20), (5.21) we obtain for all $u = (u_z, u_y) \in D$, $x = (z, y) \in D$, $w = (w_z, w_y) \in D$:

$$\begin{aligned}\nabla V(x)f(u,x,w) = &-\mu v_{\max}^2 y^2 - \left(\Phi'(w_z+\lambda) - \Phi'(z+\lambda)\right)b\left(\Phi'(w_z+\lambda) - \Phi'(z+\lambda)\right)\\&+\Phi'(\lambda+z)\left(u_y p\left(u_y, \Phi'(z+\lambda) - \Phi'(u_z+\lambda)\right) - b\left(\Phi'(z+\lambda) - \Phi'(u_z+\lambda)\right)\right)\\&-\Phi'(w_z+\lambda)\left(yp\left(y, \Phi'(w_z+\lambda) - \Phi'(z+\lambda)\right) - b\left(\Phi'(w_z+\lambda) - \Phi'(z+\lambda)\right)\right)\end{aligned} \qquad (6.33)$$

It follows from (6.33) that (4.1) holds with

$$c = \mu v_{\max}^2, \quad p = 2 \qquad (6.34)$$

$$\begin{aligned}g(u,x) := &\Phi'(\lambda+z)\left(u_y p\left(u_y, \Phi'(z+\lambda) - \Phi'(u_z+\lambda)\right) - b\left(\Phi'(z+\lambda) - \Phi'(u_z+\lambda)\right)\right)\\r(x,w) := &-\left(\Phi'(w_z+\lambda) - \Phi'(z+\lambda)\right)b\left(\Phi'(w_z+\lambda) - \Phi'(z+\lambda)\right)\\&-\Phi'(w_z+\lambda)\left(yp\left(y, \Phi'(w_z+\lambda) - \Phi'(z+\lambda)\right) - b\left(\Phi'(w_z+\lambda) - \Phi'(z+\lambda)\right)\right)\end{aligned} \qquad (6.35)$$

Using (6.35) we get for all $x = (z, y) \in D$, $w = (w_z, w_y) \in D$:

$$r(x,w) + g(x,w) = -\left(\Phi'(w_z+\lambda) - \Phi'(z+\lambda)\right)b\left(\Phi'(w_z+\lambda) - \Phi'(z+\lambda)\right) \qquad (6.36)$$

Relations (5.10), (6.36) show that (4.2) holds with $\sigma = \omega = 0$.

The fact that $S = \mathbb{R}_+ \times \mathbb{R}$ (recall (5.18)) and (6.35), (5.4) show that the following equation holds for all $x = (z, y) \in D$, $w = (w_z, w_y) \in S$:

$$r(x, w) = \Phi'(z + \lambda) b(-\Phi'(z + \lambda)) \tag{6.37}$$

Relations (5.10), (6.37) show that (4.4) holds with $\sigma = \Gamma_2 = 0$.

We next show the validity of (4.3). To this purpose, we exploit definition (5.13), the fact that $b : \mathbb{R} \to (v^* - v_{max}, v^*)$ (which implies that $0 < v^* - b(x) < v_{max}$ for all $x \in \mathbb{R}$) to get the following estimate for all $x, y \in \mathbb{R}$:

$$
\begin{aligned}
|p(y, x)| &\leq \frac{|y||v_{max} - 2v^* + 2b(x)| + \sqrt{y^2 v_{max}^2 + 4(v_{max} - v^* + b(x))(v^* - b(x))}}{2(1 + y^2)} \\
&\leq \frac{|y||v_{max} - 2v^* + 2b(x)| + \sqrt{y^2 v_{max}^2 + v_{max}^2}}{2(1 + y^2)} \leq \frac{v_{max}|y|}{2(1 + y^2)} + \frac{|y||v^* - b(x)|}{(1 + y^2)} + \frac{v_{max}}{2\sqrt{y^2 + 1}} \\
&\leq \frac{3 v_{max} |y|}{2(1 + y^2)} + \frac{v_{max}}{2\sqrt{y^2 + 1}} \leq \frac{5 v_{max}}{4}
\end{aligned}
\tag{6.38}
$$

Consequently, the conjunction of (6.35), (5.4), (6.38) and the fact that $\Omega \subseteq \mathbb{R}_+ \times \mathbb{R}$ (recall (5.18)) gives for all $u = (u_z, u_y) \in \Omega$, $x = (z, y) \in D$:

$$
\begin{aligned}
g(u, x) &= u_y \Phi'(\lambda + z) p(u_y, \Phi'(z + \lambda)) - \Phi'(\lambda + z) b(\Phi'(z + \lambda)) \\
&\leq \frac{5 v_{max}}{4} |u_y| |\Phi'(\lambda + z)| - \Phi'(\lambda + z) b(\Phi'(z + \lambda))
\end{aligned}
\tag{6.39}
$$

Define $\lim_{x \to -\infty} (b(x)) = -2\varpi$. Notice that (5.10) implies that $\varpi > 0$. Consequently, by continuity and monotonicity of $b : \mathbb{R} \to (v^* - v_{max}, v^*)$ there exists $X > 0$ with $b(-X) = -\varpi$ such that $b(x) \leq -\varpi$ for all $x \leq -X$. Define:

$$\Omega = \mathbb{R}_+ \times [-\Lambda, \Lambda] \text{ with } \Lambda := \frac{4\varpi}{5 v_{max}} \tag{6.40}$$

It follows from (6.39), (6.40) that the following inequality holds for all $u = (u_z, u_y) \in \Omega$, $x = (z, y) \in D$ with $\Phi'(z + \lambda) \leq -X$ (which implies that $b(\Phi'(z + \lambda)) \leq -\varpi$):

$$g(u, x) \leq \varpi |\Phi'(\lambda + z)| + \varpi \Phi'(\lambda + z) \leq 0 \tag{6.41}$$

Notice that (5.4), (5.5) and (6.39) implies that the following inequality holds for all $u = (u_z, u_y) \in \Omega$, $x = (z, y) \in D$ with $\Phi'(z+\lambda) \geq 0$ (which implies that $\Phi'(z+\lambda) = 0$):

$$g(u,x) \leq 0 \tag{6.42}$$

Next, we focus on the case where $0 > \Phi'(z+\lambda) > -X$. Define $\eta = \min\{b'(s) : s \in [-X, 0]\}$ ad notice that (5.10) implies that $\eta > 0$. Since $b(0) = 0$, we obtain that

$$-\Phi'(z+\lambda) b(\Phi'(z+\lambda)) \leq -(\Phi'(z+\lambda))^2 \min\{b'(s) : s \in [\Phi'(z+\lambda), 0]\} \leq -\eta (\Phi'(z+\lambda))^2$$

Therefore, we obtain from (6.39) for all $u = (u_z, u_y) \in \Omega$, $x = (z, y) \in D$ with $0 > \Phi'(z+\lambda) > -X$:

$$g(u,x) \leq \frac{5 v_{\max}}{4} |u_y| |\Phi'(\lambda+z)| - \eta (\Phi'(z+\lambda))^2 \leq \frac{25 v_{\max}^2}{64 \eta} |u_y|^2 \tag{6.43}$$

Combining (6.41), (6.42) and (6.43) we conclude that (4.3) holds with $\Gamma_1 = \frac{25 v_{\max}^2}{64 \eta}$ and $\omega = 0$.

Thus, we have shown that (4.1), (4.2), (4.3), (4.4) hold with $\Omega = \mathbb{R}_+ \times [-\Lambda, \Lambda]$, $\Lambda := \frac{4\varpi}{5 v_{\max}}$, $c = \mu v_{\max}^2$, $p = 2$, $\Gamma_1 = \frac{25 v_{\max}^2}{64 \eta}$, $\sigma = \omega = \Gamma_2 = 0$ and $r, g$ defined by (6.35), where $\varpi = \frac{1}{2} \lim_{x \to -\infty} (-b(x))$, $\eta = \min\{b'(s) : s \in [-X, 0]\}$ and $X > 0$ is defined by the equation $b(-X) = -\varpi$.

We finish the proof by applying Theorem 4 and by showing that the assumption that (2.1), (5.18), (5.19) with $\Omega = \mathbb{R}_+ \times [-\Lambda, \Lambda]$ is forward complete for every $n \in \mathbb{N}$ is not needed.

Indeed, definition (5.21) and properties (5.3), (5.4), (5.5) imply that $\{(z, y) \in (L-\lambda, +\infty) \times \mathbb{R} : V(z, y) \leq \rho, |(z, y)| \leq \rho\}$ is compact for every $\rho > 0$: the inequality $V(z, y) \leq \rho$ implies that there exists $\psi \in (L-\lambda, 0]$ such that $z \geq \psi$. Consequently, the set $\Phi_\rho = \{(z, y) \in (L-\lambda, +\infty) \times \mathbb{R} : V(z, y) \leq \rho\}$ for every $\rho > 0$ is a subset of the set $[\psi, +\infty) \times \left[-\frac{\sqrt{2\rho}}{v_{\max}}, \frac{\sqrt{2\rho}}{v_{\max}}\right]$. Therefore, definitions (5.18), (5.19), property (5.4) and the fact that $\Omega = \mathbb{R}_+ \times [-\Lambda, \Lambda]$ shows that for every compact sets $\bar{\Omega} \subseteq \Omega$, $\bar{S} \subseteq S$ and every $\rho > 0$, the sets $f(\Phi_\rho \times \Phi_\rho \times \Phi_\rho)$, $f(\bar{\Omega} \times \Phi_\rho \times \Phi_\rho)$, $f(\Phi_\rho \times \Phi_\rho \times \bar{S})$ are bounded.
The proof is complete. ◁

# 7. Concluding Remarks

The paper provided trajectory-based and Lyapunov-based sufficient conditions for string stability in homogeneous and possibly bidirectional strings. There are two research directions where the obtained results can be utilized in an instrumental way.

- First, the provided sufficient conditions can be used in conjunction with feedback design methodologies to obtain feedback laws that guarantee both asymptotic stability and string stability. Indeed, string stability in such cases guarantees robustness with respect to external inputs and *robustness with respect to the dimension of the system*. The latter characteristic cannot be guaranteed by other stability properties that take into account external inputs (like Input-Output stability, or Input-to-Output Stability).
- Since strings arise naturally when spatial discretization is applied to PDEs, it is possible that string stability can be utilized to the study of stability properties of PDEs. The key point here is that string stability is independent of the length of the string (or equivalently the step of the spatial discretization) and can allow to let the length of the string to tend to infinity. Thus, we may obtain crucial stability properties for the solution of the PDE.

Another research direction that can be studied is the derivation of necessary conditions for string stability for special classes of strings (e.g., linear strings). Novel mathematical results are needed to this purpose.

# Acknowledgments

The work of Iasson Karafyllis and Markos Papageorgiou was funded by the European Research Council under the European Union's Horizon 2020 Research and Innovation programme/ ERC Grant Agreement n. [833915], project TrafficFluid.

Iasson Karafyllis would like to thank Zhong-Ping Jiang for his friendship and his support. Zhong-Ping Jiang has been a mentor to him and their long collaboration has influenced him very much.